\documentclass[12pt]{article}
\usepackage{amsmath, amsfonts, amscd, amsthm}
\DeclareMathOperator{\lk}{lk}
\DeclareMathOperator{\cok}{cok} 
\DeclareMathOperator{\quot}{quot} 
\usepackage{epsfig}

\newtheorem{theorem}{Theorem}[section]
\newtheorem{proposition}[theorem]{Proposition}
\newtheorem{lemma}[theorem]{Lemma}
\newtheorem{diagram}[theorem]{Figure}
\newtheorem{corollary}[theorem]{Corollary}

\theoremstyle{definition}
\newtheorem{definition}[theorem]{Definition}

\newtheorem{conjecture}[theorem]{Conjecture}
\newtheorem{example}[theorem]{Example}

\newtheorem{remark}[theorem]{Remark}

\newcommand{\Q}{{\bf Q}}
\newcommand{\R}{{\bf R}}
\newcommand{\Z}{{\bf Z}} 
\newcommand{\CC}{{\bf C} [{\bf C}]}
\newcommand{\ZZ}{{\bf Z} [{\bf Z}]}
\newcommand{\Zq}{{\bf Z} [{\bf Z}/q]}
\newcommand{\Ztpm}{{\bf Z}[t^{\pm 1}]}
\newcommand{\ZqZ}{{\bf Z}[{\bf Z}/q \times {\bf Z}]}

\title{Alexander Polynomials of Equivariant Slice and Ribbon 
Knots in $S^3$}
\author{James F.~Davis\footnote{Partially supported by a grant from the
National Science Foundation} ~and Swatee Naik}
\date{}
\begin{document}
\maketitle

\begin{abstract}
This paper gives an algebraic characterization of Alexander 
polynomials of equivariant ribbon knots
  and  a factorization 
condition satisfied by Alexander polynomials of
equivariant slice knots.
\end{abstract}

\section{Introduction}
\begin{definition}
\label{slice}
A knot in $S^3$ is said to be {\em slice} if it bounds 
a slice disk, i.e., a smooth 2-disk
properly embedded in the 4-ball. It is called {\em ribbon} 
if it bounds a slice disk such that
the radius function of the 4-ball, when 
restricted to the disk, is a Morse function with critical points 
of index 0 
and 1 only. Such a disk is called a ribbon disk.
\end{definition}

Ribbon knots are alternatively defined \cite{F2,Ro} as knots which bound an 
immersed 2-disk $h : D^2 \to S^3$, so that every component of 
self-intersection is an arc $A$ so that 
$h^{-1}(A)$ consists of two arcs in $D^2$,  one of which is interior. 
Whether or not all slice knots are ribbon 
is an open question which first appeared as
Problem 25 in Fox's problem list
\cite{F2} in 1966. Also see Problem 1.33 in \cite{Ki}.

\begin{definition}
\label{periodic}
A knot $K$ in $S^3$ is said to be {\em periodic}
with period $q$ if there 
exists an order $q $ orientation preserving diffeomorphism
$f \colon S^3 \to S^3$ such that $f(K) = K$, and the fixed
point set of $f$ is a circle disjoint from $K$. 
\end{definition}

By the positive solution to the Smith Conjecture  \cite{MB}, Fix$(f)$ is 
unknotted. It follows that
the quotient space of $S^3$ under the ${\bf Z}/q$ action 
is again $S^3$. We denote the image of $K$ in the quotient space
by $\overline K$, and call it the quotient knot. We let
$B=$ Fix$(f)$ and $\overline{B}$ be its image in the quotient space, and
refer to them
as the axes. The quotient link is $\overline{B}  \cup\overline K$. 

For  example, the trefoil in its usual 
diagram is clearly seen to have 
period 3 and with an
alternate diagram we see that it 
has period 2; both quotient knots are trivial.

\begin{definition}
A period $q$ knot is called {\em $q$-equivariant slice (ribbon)}
if  the periodic diffeomorphism of $S^3$ extends to that
of $B^4$ and there is an invariant slice (ribbon) disk.
\end{definition}

A periodic  knot which bounds a 2-disk immersed
in $S^3$ with ribbon singularities, where the 2-disk is preserved by the
periodic action is equivariant ribbon (see \cite{N}). 
 An example of a 2-equivariant ribbon knot is the
Stevedore's knot $6_1$, see Figure \ref{Stevedore}. 

Clearly, equivariant
ribbon knots are 
equivariant slice. In fact,
all the knots that have been shown so far to be equivariant slice 
 are equivariant ribbon.

\begin{conjecture} \label{sliceribbon}
An 
equivariant slice knot is equivariant ribbon.
\end{conjecture}

It follows from Smith theory \cite{Smith} 
that the fixed point set of the periodic map of $B^4$
is a disk.
The fixed point disk may be knotted, as
there are non-standard finite group
actions on $B^4$  \cite{Giffen}. However, the quotient
manifold is  a homology 4-ball whose boundary is $S^3$. 
In fact, one can show that the quotient manifold is a simply-connected
homology 
4-ball whose boundary is $S^3$, and hence, by topological surgery, that it is
homeomorphic to a 4-ball.

It was shown in Corollary 3.3 of \cite{N} that 
the linking number of the equivariant slice knot
with its axis is always 1. Further obstructions to a slice,
 periodic
knot  being equivariant slice were
obtained in \cite{CK, KCS, N} in terms of 
Seifert matrices, metabolizers of Seifert forms, Casson-Gordon 
invariants,
and surgery on the quotient knot. 
In \cite{CK, KCS} examples were given of knots which are slice, periodic,
and have linking number one with the axis, but cannot be 
equivariant 
slice.  (See also Example \ref{10_{123}}.)

The goal of this paper is to characterize Alexander polynomials of
equivariant
ribbon knots and to give a necessary condition that an Alexander polynomial
of
an equivariant slice knot must satisfy.  Our results place restrictions
on which knots can be equivariant slice/ribbon and construct new examples
of
equivariant ribbon knots. 

Alexander polynomials of knots, of slice knots, and of periodic knots with
linking number 1 have been characterized, and before we state our results
we
review this earlier work.  We first make a convention: all polynomials in the
Laurent polynomial ring $\Z[t,t^{-1}]$ are only defined up to  multiplication
by $\pm
t^i$.  Similarly, a two variable polynomial in $\Z[\Z \times \Z] =
\Z[g,g^{-1},t, t^{-1}]$ or in $\Z[\Z/q \times \Z] = \Z[g,t,t^{-1}]/\langle g^q
\rangle$ is only
defined  up to  multiplication by $\pm  g^i t^j$.  All equalities involving
polynomials will only hold up to such multiplies.

\begin{definition} A polynomial $\Delta(t) \in \Z[t,t^{-1}]$
is an {\em (abstract) Alexander polynomial} if
\begin{enumerate}
\item $\Delta(1) = 1$
\item $\Delta(t) = \Delta(t^{-1})$
\end{enumerate} 
\end{definition}

\begin{proposition}[Seifert, see \cite{Levine}]  An Alexander polynomial of a
knot is an
abstract Alexander
polynomial.  If $\Delta(t)$ is an abstract Alexander polynomial, then there is
a
knot whose Alexander polynomial is $\Delta(t)$.
\end{proposition}

The corresponding result for slice and ribbon knots is stated next.
Proposition \ref{slicep} below combines the result of Fox \cite{F, FM}
regarding factorization of 
the Alexander polynomial of a slice knot
together with a result of 
Terasaka \cite{Te} that, given 
a polynomial satisfying the Fox factorization condition, there exists a slice,
in fact, a
ribbon knot with the polynomial as its Alexander polynomial.

\begin{proposition}
\label{slicep}
If $\Delta(t)$ is the Alexander polynomial of a slice knot, then there is a
polynomial $p(t) \in {\bf Z}[t,t^{-1}]$
so that $\Delta(t) = p(t)p(t^{-1})$. 
Conversely, given a polynomial $p(t) \in {\bf Z}[t,t^{-1}]$ with $p(1) = 1$,
there exists a
ribbon knot whose Alexander polynomial is $p(t)p(t^{-1})$.
\end{proposition}

Alexander polynomials of periodic knots which have linking number one
with the axis have been characterized.
Proposition \ref{Murasugi} below  states the result of
K. Murasugi \cite{M1} regarding the factorization of
the Alexander polynomial of a
periodic knot, together with the existence result from \cite{DL}.

\begin{proposition} \label{Murasugi}
Let $K$ be a period $q$ knot with linking number one.  Let $\Delta_K(t)$ and
$ \Delta_{\overline K}(t)$
be the 
Alexander polynomials of $K$ and the quotient knot $\overline K$
respectively.
Then there exists  a polynomial
$\Delta_{\Z /q} (g,t) \in \ZqZ$ such that
\begin{enumerate}
\item $\Delta_{\overline K}$ divides $\Delta_K$,
\item $\Delta_{K} (t)/  \Delta_{\overline K} (t)
= \prod^{q-1}_{i=1} \Delta_{\Z /q}(\zeta^i,t)$,
where $\zeta$ is a primitive $q$-th root of unity,
\item $\Delta_{\Z /q} (g,t) = \Delta_{\Z /q} (g^{-1}, t^{-1}),$
\item $\Delta_{\Z /q} (g,1)=1,$
\item  $\Delta_{\Z /q} (1,t)= \Delta_{\overline K}(t) .$ 
\end{enumerate}

Conversely, given a polynomial $\Delta_{\Z /q} (g,t) \in \ZqZ$ satisfying 
conditions (3)~and (4)~above, 
there is a period $q$ knot with linking number one, with the Alexander
polynomial
 of the quotient knot given by (5), and the Alexander polynomial of the 
knot given by (2). 

\end{proposition}

For a periodic knot, the polynomial $\Delta_{\Z /q}(g,t) $ 
is the image  of the 2-variable Alexander
polynomial of the quotient link under the map $\Z[\Z \times \Z] \to \ZqZ$.
See \cite{M1,DL}.
We will call $\Delta_{\Z /q}(g,t) $ the {\it Murasugi polynomial} of the
periodic
knot.

A consequence of the fact that a ribbon disk has only 
critical points of index 0 and 1 is that
the exterior $X$ of a ribbon disk
 is $B^4$ with 1-handles and 
2-handles attached. Clearly, $X$ has the homotopy type of a 2-complex;
a homological consequence is that the infinite cyclic cover 
$\widetilde X$ of $X$ has trivial $H_2$. (Since $X$ is a 2-complex, $H_2
(\widetilde
X)$ is torsion-free
over
$\Z [t, t^{-1}]$, whereas the Wang sequence of the infinite cyclic cover 
shows that $H_2 (\widetilde
X)$ must be torsion  over $\Z [t, t^{-1}]$.)
This need not be the case if $X$ was merely the exterior of a slice 
disk.  This homological difference between slice disks and ribbon disks is
reflected in our
equivariant results given below.

In Section \ref{Application to equivariant slice/ribbon knots}, 
we prove the following condition on Alexander polynomials
of equivariant slice knots
using an argument based on Reidemeister torsion.

\begin{theorem}
\label{eqslicethm}
Let $K$ be a $q$-equivariant slice knot with  Murasugi polynomial $\Delta_{\Z
/q}(g,t)$.
Then there exist nonzero polynomials $a(g,t)$ and $b(g,t)$
in $\ZqZ$ such that
\begin{enumerate} 
\item $\Delta_{\Z /q} (g,t)b(g,t)b(g^{-1},t^{-1}) = a(g,t)a(g^{-1},t^{-1})
 $ 
\item $a(g,1) = 1 = b(g,1)$.
\end{enumerate}
 \end{theorem}

The fact that $H_2$ of the infinite cyclic cover of the exterior
of a ribbon disk is trivial enables us to show that if $K$ is
equivariant ribbon, then $b(g,t) = 1$.
We show that this characterizes 
Alexander polynomials of equivariant ribbon knots. 
Our result is stated below. The proof
uses a handle-theoretic construction which is described in Section
\ref{Construction of
equivariant ribbon knots}.

\begin{theorem}
\label{eqribbonthm}
Let $K$ be a $q$-equivariant ribbon knot with Murasugi polynomial
$\Delta_{\Z /q}(g,t)$. Then there
exists a polynomial 
$a(g,t)$ 
in $\ZqZ$ such that
\begin{enumerate}
\item 
$\Delta_{\Z /q} (g,t) = a(g,t)a(g^{-1},t^{-1})$, and 
\item $a(g,1) = 1 .$ 
\end{enumerate}

Conversely, given a polynomial $a(g,t) \in \ZqZ$ with $a(g,1) = 1$, 
there is a $q$-equivariant ribbon knot whose Murasugi polynomial is  
given by (1), and hence whose Alexander polynomial is given by
$\prod^{q-1}_{i=0} \Delta_{\Z /q}(\zeta^i,t)$.
\end{theorem}

It follows from this result that 
in the examples of \cite{CK, KCS}, the knots have Alexander polynomials
of equivariant ribbon knots.

In general, algebraic considerations do not seem to tell
us that $b(g,t)$ has to be 1; in fact, we have an example (see Proposition 
\ref{2prop})  
 of a 2-equivariant knot whose Murasugi polynomial satisfies the conditions 
\ref{eqslicethm}(1),(2), but for which 
$b(g,t)$ cannot be chosen to be 1.  Nonetheless,
we state below a conjecture weaker than Conjecture \ref{sliceribbon}:

\begin{conjecture} 
Alexander polynomials of equivariant slice knots satisfy the 
factorization condition in \ref{eqslicethm} with
$b(g,t) = 1$. \end{conjecture}

\begin{example}
The Stevedore's knot $6_1$ pictured below is 2-equivariant ribbon.  
\begin{diagram}
\label{Stevedore}
$$\epsfxsize 3in
 \epsfysize 2.5in 
\epsffile{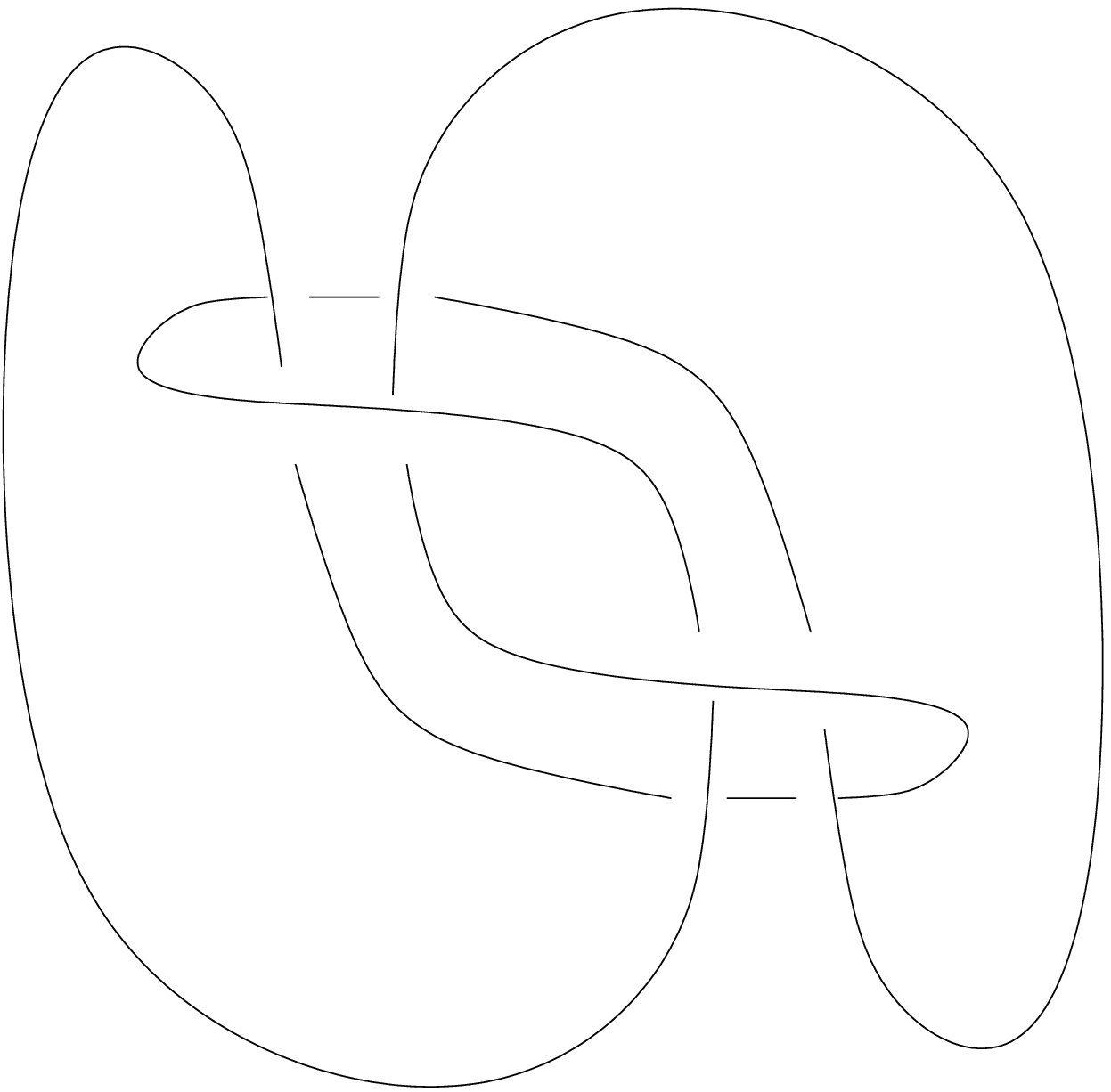}$$
\end{diagram}
Here are the various polynomial invariants for Stevedore's knot.

\begin{align*}
\Delta_K(t) &=  -2t + 5 -2t^{-1} \\
p(t) &= -2t+1\\
 \Delta_{\overline K}(t) &= 1\\
\Delta_{\Z/2}(g,t) &= (g-1)t + 3-2g + (g^{-1}-1)t^{-1} \\
a(g,t) &=  (g-1)t +1  \\
\end{align*}
\end{example}

The Alexander polynomial of $10_{123}$ satisfies the conditions
for being slice and period 2.
Our results show that it 
is not a 2-equivariant slice polynomial. 
(See Example \ref{10_{123}}.) Note that this knot
is slice \cite{Conway} and it is listed as period 2 in
\cite{BZ}, but 
in an email 
communication with the second author 
Jeff Weeks has told us that according to Snappea,
$10_{123}$ does not have period 2. 

In \cite{N} a list of periodic, slice knots, 
which are not equivariant slice, was given.
There were some cases of knots with less than or equal to 10 
crossings that the techniques in \cite{N} were unable
to settle. 
It follows from our results that
the Alexander polynomial of 
each  of these knots (except $10_{123}$) is that
of an equivariant ribbon knot.  This follows from 
 Corollary \ref{modq}  which states that
if a slice polynomial is 1 mod $q$, for a prime $q$,
then there exists a $q$-equivariant ribbon knot with that as its
polynomial. 

In Section \ref{Construction of equivariant ribbon knots},  
given a polynomial which satisfies the factorization conditions
in Theorem \ref{eqribbonthm}, we describe a construction of 
a $q$-equivariant ribbon knot with that as its polynomial.
A special case of our result is the known fact that 
given a polynomial of the form $p(t)p(t^{-1})$, there exists
a ribbon knot with that as its polynomial. Our proof is handle-theoretic
and avoids the long computations of \cite{Te}.
(Handle-theoretic constructions of slice knots was given earlier in \cite{AS}.)
In Section \ref{Application to equivariant slice/ribbon knots} we 
give an argument based on Reidemeister torsion 
to obtain necessary factorization conditions for Alexander polynomials
of equivariant slice knots and those of equivariant ribbon knots.

\section{Applications} \label{Applications}

We obtain below a realization result as a consequence of Theorem 
\ref{eqribbonthm}. (Compare Corollary 1.2 of \cite{DL}.)

\begin{corollary} \label{modq}

If an Alexander polynomial $\Delta (t)$ factors as
$p(t) p(t^{-1})$ with $p(t) \equiv 1\ {\rm mod\ } q $, 
 then there exists a $q$-equivariant ribbon knot
$K$ such that $\Delta (t)$ is the polynomial of $K$ as well as 
that of the quotient knot.
\end{corollary}

\begin{proof} 
Write $p(t) = q \cdot h(t) + 1$. 
Now let $a(g,t) = (1 + g + \cdots + g^{q-1}) \cdot h(t) +1,$
and apply Theorem \ref{eqribbonthm}.
\end{proof}

Note that if $q$ is prime, $(\Z/q)[t]$ is a unique factorization domain, so
that if $\Delta(t) = p(t)p(t^{-1})$ and if
$\Delta(t)
\equiv 1\ {\rm mod\ } q $, then it is automatic that $p(t) \equiv 1\ {\rm
mod\ } q $.

\begin{example}
\label{ex1}
The knots $8_8, \ 9_{46},\ 10_{22},$ and
$ 10_{35}$ are all slice and have period 2 with $\lambda = 1$. 
For each of these knots  
$\Delta (t) \equiv 1\ {\rm mod\ } 2.$ Hence by Corollary
\ref{modq} it follows that each of their Alexander polynomials is the
Alexander polynomial of a 
2-equivariant ribbon knot.
It was shown in \cite{KCS}, using a characterization of the Seifert
matrix, that $8_8$ is not 2-equivariant slice. It is not yet 
known whether or not any of the others are 2-equivariant slice.
\end{example}

\begin{example}
In \cite{CK} an example is given of a knot which is period 2, slice,
but not 2-equivariant slice. The polynomial of this knot is
$-2t +5 -2t^{-1}$ which is $1$ mod $2$, and hence there exists a 
2-equivariant ribbon knot with this polynomial.  (See also Figure
\ref{Stevedore}.)
\end{example}

The following is an immediate 
consequence of Theorem \ref{eqslicethm}.

\begin{corollary} \label{conjfac}
Let $K$ be a $q$-equivariant slice knot with $q$ a prime. 
Let $\zeta$ be a primitive $q$-th root of unity.
Then $\Delta_K(t) /  \Delta_{\overline K} (t)$ factors into $q-1$ conjugate
factors
over ${\bf Z} [\zeta]$ each of which is of the form 
$\left( f_1(t) f_1(t^{-1}) \right) / \left( f_2(t) f_2(t^{-1}) \right)$, 
where $f_1(t), f_2(t) \in {\bf Z} [\zeta ][t]$.
\end{corollary}

In the cases where $Z[\zeta]$ is a unique factorization domain, this can 
be simplified using the following lemma.

\begin{lemma} \label{factor} Let $R$ be a unique factorization domain.  If
$$f(t)b(t)b(t^{-1})
= a(t)a(t^{-1}) \in R[t,t^{-1}],$$ where $b(t)$ is nonzero, then   there is a
$c(t)\in R[t,t^{-1}]$ so
that 
$$f(t) = c(t)c(t^{-1}) \in R[t,t^{-1}].$$
\end{lemma}

\begin{proof}  Let $b(t) = p_1(t)^{e_1} \cdots p_r(t)^{e_r}$ be the prime
factorization of $b(t)$. Let $f_i$ be the largest integer so that $f_i \leq
e_i$ and  $p_i(t)^{f_i}$ divides $a(t)$. 
Let 
$$
c(t) = \frac{a(t)}{\prod_{i=1}^r p_i(t)^{f_i} p_i(t^{-1})^{e_i-f_i}}.
$$
\end{proof}

In the case of period 2, the factorization mentioned in \ref{conjfac}
is over the integers, hence the algebraic conditions on Alexander
polynomials of 2-equivariant slice/ribbon knots can be simplified. 
We now embark on a fairly complete discussion of this case.

First we introduce some more notation.

\begin{definition} 
Let $\Delta(t)$ and
 $\Delta^{\quot} (t)$ be abstract Alexander polynomials so that
$\Delta^{\quot} (t)~|~\Delta(t)$.   Let $q$ be a natural number.
Then $(\Delta,\Delta^{\quot})$ satisfies the {\em $q$-equivariant slice
conditions} if there exists polynomials
$\Delta_{\Z_q}(g,t), a(g,t), b(g,t) \in \ZqZ$ so that  
\begin{enumerate}
\item $\Delta_{\Z /q} (g,t)b(g,t)b(g^{-1},t^{-1}) = a(g,t)a(g^{-1},t^{-1})$
\item $a(g,1) = 1 = b(g,1)$
\item $\Delta(t) = \prod^{q-1}_{i=0} \Delta_{\Z /q}(\zeta^i,t),
\text{ where $\zeta$ is a primitive $q$-th root of unity}$
\item $\Delta^{\quot}(t) = \Delta_{\Z /q}(1,t)$
\end{enumerate}
and $(\Delta,\Delta^{\quot})$ satisfies the {\em $q$-equivariant ribbon
conditions} if there exists a polynomial $a(g,t) \in
\ZqZ$ so that 
\begin{enumerate}
\item $a(g,1) = 1$
\item $\Delta(t) = \prod^{q-1}_{i=0} a(\zeta^i,t)a(\zeta^{-i},t^{-1}),
\text{ where $\zeta$ is a primitive $q$-th root of unity}$
\item $\Delta^{\quot}(t) = a(1,t)a(1,t^{-1})$.
\end{enumerate}
\end{definition}

We can now rephrase our main theorems, \ref{eqslicethm} and
\ref{eqribbonthm} as saying.
\begin{corollary}
\begin{enumerate}
\item If $K$ is a $q$-equivariant slice knot, then
$(\Delta_K,\Delta_{\overline K})$ satisfy the 
$q$-equivariant slice conditions.
\item  If $K$ is a $q$-equivariant ribbon knot, then
$(\Delta_K,\Delta_{\overline K})$ satisfy the 
$q$-equivariant ribbon conditions. Furthermore, if $(\Delta,\Delta^{\quot})$
satisfy the $q$-equivariant ribbon
conditions, there is a $q$-equivariant ribbon knot $K$, so that\linebreak
$(\Delta,\Delta^{\quot}) =
(\Delta_K,\Delta_{\overline K})$.
\end{enumerate}
\end{corollary}

\begin{proposition} \label{2prop}
\begin{enumerate}
\item  $(\Delta,\Delta^{\quot})$ satisfies the 2-equivariant slice conditions
if and only if there are polynomials
$p(t), q(t)$ so that
\begin{enumerate}
\item $\Delta(t) = p(t)p(t^{-1})$
\item  $\Delta^{\quot}(t) = q(t)q(t^{-1})$
\item $q(t)~|~p(t)$
\item $\Delta(t) \equiv ( \Delta^{\quot}(t) )^2 \pmod{2}$.
\end{enumerate}
\item  $(\Delta,\Delta^{\quot})$ satisfies the 2-equivariant ribbon
conditions if and only if there are polynomials
$p(t), q(t)$ so that
\begin{enumerate}
\item $\Delta(t) = p(t)p(t^{-1})$
\item  $\Delta^{\quot}(t) = q(t)q(t^{-1})$
\item $q(t)~|~p(t)$
\item $p(t) \equiv q(t)^2 \pmod{2}$
\end{enumerate}
\item There is a pair $(\Delta,\Delta^{\quot})$ satisfying the 2-equivariant
slice conditions which does not satisfy
the 2-equivariant ribbon conditions.  In fact, for all divisors $D$ of $\Delta$,
$(\Delta, D)$ does not satisfy
the 2-equivariant ribbon conditions.
\end{enumerate}
\end{proposition}

\begin{proof}  For the purposes of this proof, we shall adopt the following
notation.  For $f(t) \in \Z[\Z]$ and 
$h(g,t) \in \Z[\Z/2 \times \Z]$, let $\overline f(t) = f(t^{-1})$ and
$\overline h(g,t) = h(g^{-1},t^{-1})$.  Also, $c
\equiv d$ will mean $c \equiv d \pmod{2}$.

A polynomial $h(g,t) \in \Z[\Z/2 \times \Z]$ determines polynomials $h_+(t)
= h(1,t)$ and $h_-(t) = h(-1,t)$ so that
$h_+(t) \equiv h_-(t)$.  Conversely, polynomials $h_+, h_- \in \Z[\Z]$ with
$h_+ \equiv h_-$ determine a polynomial 
$$h(g,t) =
\frac{1+g}{2}h_+(t) + \frac{1-g}{2}h_-(t).
$$

\noindent (1):  Suppose  $(\Delta,\Delta^{\quot})$ satisfies the 
2-equivariant slice conditions, and thus $\Delta_{\Z
/2}(g,t), a(g,t), b(g,t)$ are given.  Applying Lemma \ref{factor} to the
equation $(\Delta_{{\Z/2}-}
)(b_-)(\overline{b_-}) = (a_-)(\overline{a_-})$, we see $\Delta_{{\Z/2}-} =
q\overline q$ for some $q$.  Likewise 
$\Delta_{{\Z/2}+} = r\overline r$ for some $r$.  Then 
desired conditions are satisfied by setting $p = rq$.

Suppose now $(\Delta,\Delta^{\quot})$ satisfies  conditions 1(a)-1(d) above
and thus $p(t)$ and $q(t)$ are given.
Let $r = p/q$.  Then condition 1(d) implies that $r\overline r \equiv
q\overline q$.  Using that $(\Z/2)[\Z]$ is a
unique factorization domain and proceeding by induction on the degree of
$r$, one can prove that there are polynomials
$c(t), d(t)
\in
\Z[\Z]$ so that
$c(1) = 1 = d(1)$,
$r \equiv cd$, and $q \equiv c \overline d$.  Thus $r \overline c \overline d
\equiv q \overline c d$.  Letting $b(g,t)
= c(t)d(t)$ and 
$$
a(g,t) = \frac{1+g}{2}q(t) \overline c(t) d(t)  + \frac{1-g}{2}r(t) \overline
c(t)\overline d(t).
$$

\noindent (2):  Suppose  $(\Delta,\Delta^{\quot})$ satisfies the 
2-equivariant ribbon conditions, and thus $a(g,t)$ is
given.  Then let $q= a_+$ and $p = a_+ a_-$.

Suppose now $(\Delta,\Delta^{\quot})$ satisfies  conditions 2(a)-2(d) above
and thus $p(t)$ and $q(t)$ are given.  Let
$r = p/q$.  Then condition 2(d) implies that $r \equiv q$.  Let $a(g,t) =
((1+g)/2)q(t)  +
((1-g)/2)r(t)$.

\medskip

\noindent (3):  Note that to specify $a$, $b$, $\Delta$, and $\Delta^{\quot}$ we need
polynomials $a_+, a_-, b_+,b_- \in
\Z[t,t^{-1}]$ so that 
so that 
\begin{align*}
b_+\bar b_+ &~|~a_+ \bar a_+\\ 
b_-\bar b_- &~|~a_- \bar a_-\\ 
a_+ &\equiv a_- \\
b_+ &\equiv b_- \\
a_+(1) & = a_-(1) = b_+(1)  = b_-(1) =1\\
\end{align*}
Then let $\Delta_+ = (a_+\bar a_+)/(b_+\bar b_+)$ and $\Delta_- = (a_-\bar a_-)/(b_-\bar b_-)$. 
Then $\{a_+,a_-\}$, $\{b_+,b_-\}$, and $\{\Delta_+,\Delta_-\}$ determine $a$, $b$ and $\Delta$ 
respectively, and we set $\Delta^{\quot} = \Delta_+$.

Now suppose that $f \in \Z[t,t^{-1}]$ is irreducible, that $f(1) = 1$ and that
$f \not \equiv
\bar f$ (e.g.
$f(t) = t^3 - t^2 + 1$).  Let $\alpha$ be a irreducible polynomial so that
$\alpha(1) = 1$ and 
$\alpha \equiv \bar f f$  (e.g. $\alpha(t) = t^6 -t^5 -t^4 +3t^3 - t^2 - 3t +
3)$.  Let
$\beta$ be a irreducible polynomial so that $\alpha(1) = 1$ and $\beta \equiv 
f f$  (e.g. $\beta(t) = t^6 -2t^5 +t^4
+2t^3 - 4t^2 + 3)$.  Let
\begin{align*}
b_+ &= f \\
b_- &= f \\
a_+ &= f \alpha\\
a_- & = \bar f \beta
\end{align*}
If  $\Delta = c \bar c$, then $c_+ \bar c_+ = (a_+ \bar a_+)/(b_+ \bar
b_+) = \alpha \bar \alpha$, so $c_+ = \alpha$ or $c_+ =
\bar \alpha$.    Likewise $c_- = \beta$ or $c_- =
\bar \beta$. In all cases $c_+ \not \equiv c_-$ which contradicts the
existence of $c$.

\end{proof}

\begin{example}
\label{10_{123}}
Let $\Delta_K (t) = (t^4-3t^3+3t^2-3t+1)^2.$ 
(This is the polynomial of the knot
$10_{123}$.) 
This is the Alexander polynomial of 
a slice knot and that of a period 2 knot.
Since  $\Delta_{\overline K}( t)$
and $\Delta_K(t) / 
\Delta_{\overline K}(t)$ are congruent mod 2, they are both equal to
$ t^4-3t^3+3t^2-3t+1$, which does not factor as $q(t) q(t^{-1})$
over {\bf Z}.  It follows from Proposition \ref{2prop} that
this is not an Alexander polynomial of a 2-equivariant slice
knot. 
\end{example}

Referring to the list of knot concordances given in \cite{Conway},
examples \ref{ex1} and \ref{10_{123}}
above, together with examples 2.5 and 3.4 in
\cite{N} list all the candidates for equivariant slice knots from
the tables of knots with 10 or fewer crossings.

\section{Construction of equivariant ribbon knots} \label{Construction of
equivariant ribbon knots}
In this section, we will prove the realization part of Theorem
\ref{eqribbonthm}.
For any $a(g,t) \in \ZqZ$ with $a(g,1) = 1$, 
we construct  
a $q$-equivariant ribbon knot $K$
with Murasugi polynomial $a(g,t)a(g^{-1},t^{-1})$.  We do this by
constructing the 
quotient link $\overline B \cup \overline K$.

\begin{remark} 
Let $W$ be a 4-manifold with boundary.  By {\em attaching a $j$-handle to
$W$} (see \cite[IV.6]{Ko}),
one means constructing the manifold
$$
W' = W \cup_\phi D^j \times D^{4-j} 
$$
where $\phi : S^{j-1} \times D^{4-j} \hookrightarrow \partial W$ is an
embedding.  A smooth structure on 
$W$ leads, by ``rounding corners'', to a smooth structure on $W'$, which is
invariant under isotopy of $\phi$.  
To obtain $\partial W'$  from $\partial W$   one proceeds by the process of
surgery, i.e.
$$
\partial W' \cong (\partial W - \phi ( S^{j-1} \times \text{int }   D^{4-j}))~
\cup_{S^{j-1} \times
S^{3-j}}~ D^j \times S^{3-j}
$$

If one attaches a 1-handle and then a 2-handle
$$
W'' = (W \cup_{\phi_1} D^1 \times D^{3}) \cup_{\phi_2} D^2 \times D^2,
$$
so that the belt sphere $\phi_1(1 \times S^2)$ of the 1-handle intersects
the attaching sphere $\phi_2(S^1 \times 0)$ of
the 2-handle exactly once, transversely in $\partial W'$, then $W$ and $W''$
are diffeomorphic (see
\cite[IV.7]{Ko}).
     
If $K$ is a slice knot in $S^3$ (e.g. the unknot), and if one attaches a
cancelling pair of handles as above to $B^4$
whose attaching maps $\phi_1, \phi_2$ have image disjoint from $K$, then
$K \subset \partial W'' \cong S^3$
 is still a slice knot, however, possibly a different knot from the original knot
$K
\subset S^3$.
\end{remark}

\begin{remark} \label{construction remark}  Let $D_{\overline{B}}$ and
$D_{\overline{K}}$ be a pair
of linear, orthogonal 2-disks in $B^4$.  Then
$$
\overline{B} \cup \overline{K} =
(D_{\overline{B}} \cup D_{\overline{K}}) \cap S^3 \subset S^3
$$
is the Hopf link.  Suppose one attaches a cancelling pair of handles
$\phi_1,\phi_2$ as above to $B^4$, with
$\overline{K}$ disjoint from the handles.  Suppose, furthermore, that
$\overline{B}$ bounds a 2-disk
in $S^3$, disjoint from the handles.   Let $B$, $K$, $D_B$,
and $D_K$ be the inverse images of $\overline{B}, \overline{K},
D_{\overline{B}}, D_{\overline{K}}$ under the
$q$-fold branched cover $\pi: \widehat { W''} \to \ W''$
branched over the unknotted 2-disk
$D_{\overline{B}}$.  This gives us a $q$-equivariant slice (in fact, ribbon)
knot $K  \subset \widehat{\partial
W''} \cong S^3$. 
\end{remark}

We will use Akbulut's modification of Kirby diagrams (see
\cite[Chapter 5]{GS}).  A pair of cancelling handles is indicated by a 
2-component link with one
component marked with a dot and bounding a 2-disk (usually not explicitly
drawn), which intersects
the other component transversely once. The other component is equipped
with a integer
$k$.  The component with the dot corresponds to adding a 1-handle; the
second component corresponds to adding the
2-handle, with the $k$ indicating the framing.   An arc $A$ of the undotted component
intersecting the
spanning disk for the dotted
component indicates that a piece of the attaching
map of the 2-handle goes across the boundary of the 
1-handle, i.e., $A \times D^2 \subset D^1 \times S^2 \subset
D^1 \times D^3$.

We illustrate  our construction of the quotient link
through a simple example.

\begin{example}
\label{exequiv} 
Let $a (g,t) = 1 -g +gt$. 

Figure \ref{1-g+gt} indicates the link 
$\overline{B} \cup \overline K \subset S^3$
and a pair of a 1-handle and a 2-handle
attached to the 4-ball along $S^3$.

\begin{diagram}
\label{1-g+gt}
$$
 \epsfysize 2in
\epsffile{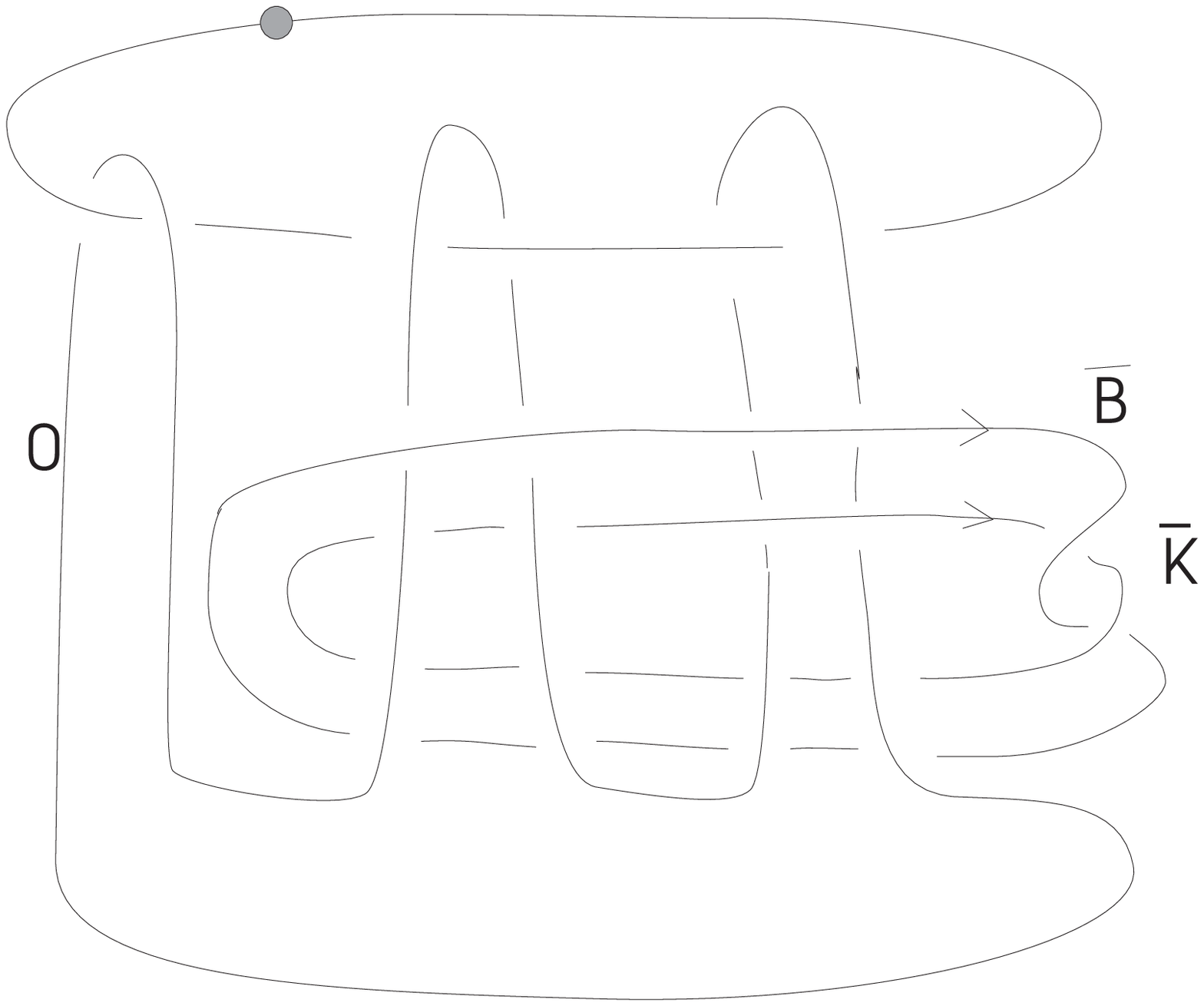}$$
\end{diagram}

\begin{remark} \label{conditions}  Let $S_1$ be the circle labelled with 0, and
let $S_2$ be the circle labelled with the dot.   
\begin{enumerate}
\item   There is a disk $D$ with boundary $S_1$, so that $S_2$ intersects  $D$   exactly
once, transversely.
\item  $\overline{B} \cup \overline K$ is the Hopf link and is disjoint from
$S_1 \cup S_2$.
\item After erasing $\overline K$, $\overline B$ bounds a 2-disk in the
complement of the rest
of the link and $D$.
\item $S_1$ is contained in a ball disjoint from  $\overline{B} \cup \overline
K$.  $S_2$ is contained in a
ball disjoint from  $\overline{B} \cup \overline K$.
\end{enumerate}
\end{remark}

As explained in Remark \ref{construction remark}, (1), (2), and (3) imply
that the $q$-fold branched cover,
branched over $\overline B$, gives a $q$-equivariant slice knot.  It is
equivariant ribbon, since the
complement of the slice disk is constructed from $B^4$ by using only 1-and
2-handles.  Alternatively,
one can see that the knot is ribbon by examining Figure \ref{ribbon2} below. 
Item (4) is convenient
for doing computations by lifting handles to the infinite cyclic cover of the
knot complement.

What remains in this example is the computation of the Murasugi polynomial. 
This can be
accomplished in three ways: by lifting handles (this is what we do next), by
using 
Poincar\'e-Lefschetz duality and Reidemeister torsion (see Section
\ref{Application to equivariant
slice/ribbon knots}), or by direct analysis of Figure \ref{ribbon2} below.

Let $\pi : X \to \overline X$ be the $q$-fold branched cover of $\overline X =
S^3 - \overline K$, branched
over $\overline B$.  Thus $X = S^3 - K$.  Let $\widetilde X \to X$ be the
infinite cyclic cover of $X$. 
Since the action of $\Z/q$ on $X$ has a nonempty fixed set, the action lifts
to $\widetilde X$. Then 
$\Z/q \times \Z$ acts on $\widetilde X$ with quotient $\overline X$.  Let $q:
\widetilde X \to \overline X$ be the
quotient map.  Note that before the handles are added, $\widetilde X = \R
\times \R^2$.

Let $S_1$ and $S_2$ be the closed curves indicating the 1- and 2-handles
respectively.  Then
$q^{-1}(S_1)$ and  $q^{-1}(S_2)$ are a disjoint collections of closed curves,
freely and transitively 
permuted by $\Z/q \times \Z$.  

\begin{definition} 
The {\em equivariant linking  of $S_1$ and $S_2$} 
$$
\lk_{\Z/q \times \Z}(S_1, S_2) = \sum_{i,j\in \Z} \lk(\widetilde S_1,
g^it^j\widetilde S_2) g^it^j \in \ZqZ
$$
is defined by choosing a component $\widetilde S_i$ of  $q^{-1}( S_i)$ and an
orientation for that component. 
The equivariant linking  is well-defined up to multiples of $\pm g^it^j$.
\end{definition}

If $-: \ZqZ \to \ZqZ$ is the standard involution
$$
\overline{\sum a_{i,j}g^it^j} = \sum a_{i,j}g^{-i}t^{-j},
$$
then the equivariant linking form satisfies the symmetry condition
$$
\lk_{\Z/q \times \Z}(S_2, S_1) = \overline{\lk_{\Z/q \times \Z}(S_1, S_2)}.
$$
The equivariant linking  can be computed from a diagram such as Figure
\ref{1-g+gt} as follows. 
First choose orientations of $S_1$ and $S_2$ and then choose a point $P_0
\in S_1$ corresponding to an
overcrossing of $S_1$ (over $S_2$).
Then
$$
\lk_{\Z/q \times \Z}(S_1, S_2) = \sum_{P} \varepsilon_P
g^{\lk(C_P,\overline B)}t^{\lk(C_P,\overline K)}
$$
where the sum is over all points where $S_1$ crosses over $S_2$.  Here
$\varepsilon_P = \pm 1$ where the
sign corresponds to whether the crossing is positive or negative.  Here
$C_P$ is the curve which starts at $P_0$,
travels along $S_1$ following the orientation until it reaches $P$ (which is at
an overcrossing), switches down to
$S_2$, follows
$S_2$ against the orientation until reaching the point under $P_0$, and then
jumps up to $P_0$.    The reader
should now study Figure \ref{1-g+gt} to check that the equivariant linking  is
$1-g+gt$.  The proof that
the equivariant linking can be computed in this way is similar to the
corresponding result for the
classical linking number \cite{Ro}. 

\begin{theorem} \label{Murasugicompute} Suppose one has a 4-component
link satisfying the conditions of Remark \ref{conditions}. 
Let $a(g,t) = \lk_{\Z/q \times \Z}(S_1, S_2)$.  Then the Murasugi polynomial
of the corresponding
$q$-equivariant ribbon knot is $a(g,t)a(g^{-1},t^{-1})$.

\begin{proof}  The boundary of the 4-manifold given by adding the 1-handle
along $S_1$ and the 2-handle
along $S_2$ is the same as the 3-manifold given by doing 0-framed surgeries
along the curves $S_1$ and
$S_2$ (see \cite[p. 171--172]{GS}).  The advantage of doing this is that the
first homology of the
3-manifold is presented by the linking matrix of the surgery curves (see
\cite[p. 165]{GS}).  This all works
equivariantly, and one sees that if $\partial W''$ denotes the result of the
surgeries and $\widetilde X''$
denotes the infinite cyclic cover of the  complement of the equivariant knot,
then 
$$
H_1(\widetilde X'') \cong \cok (
\begin{pmatrix}
0 & a(g,t) \\
a(g^{-1},t^{-1}) & *
\end{pmatrix} : 
 \ZqZ^2 \to \ZqZ^2).
$$
The first homology remains unchanged after removing the inverse image of
the branch set.  The Murasugi
polynomial is thus the determinant of the above matrix.
\end{proof}

\end{theorem}

We now return to our example (Figure \ref{1-g+gt}), and draw a picture of
the resulting ribbon knot. 
This is not, strictly speaking, logically necessary, but nonetheless
psychologically satisfying.   

By isotopy
the picture can be changed to the following.

\begin{diagram}
\label{ribbon1}
$$
 \epsfysize 2.5in
\epsffile{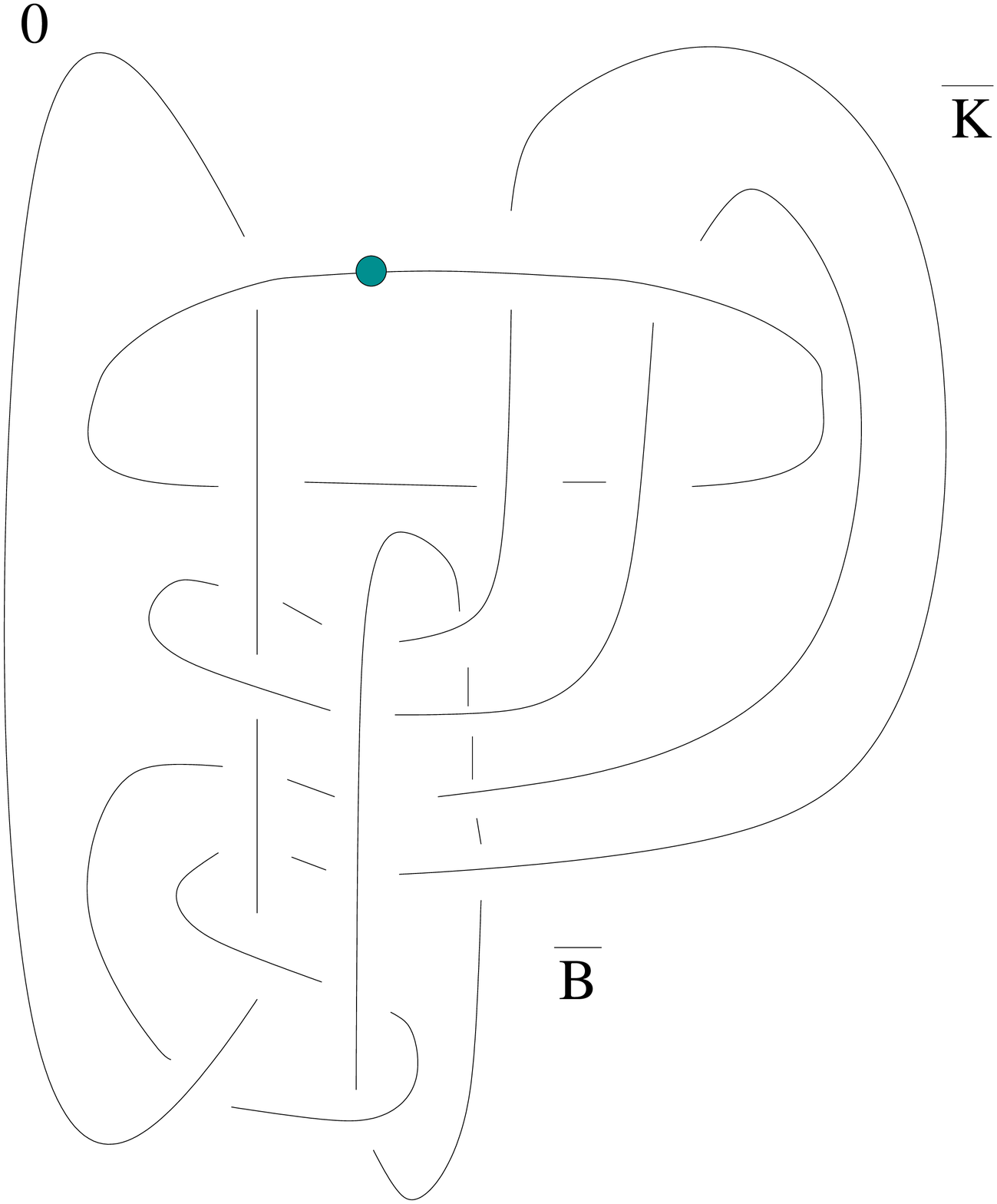}$$

\end{diagram}

In order to see the image of $\overline K$ and $\overline B$ in the
boundary of the 
new $B^4$, eliminate the intersections of $\overline K$ with the horizontal
disk $D$ by sliding $\overline K$ over the
2-handle and then erase the canceling pair of the 1-handle and the 2-handle.
We see the following.

\begin{diagram}
\label{ribbon2}
$$
 \epsfysize 2in
\epsffile{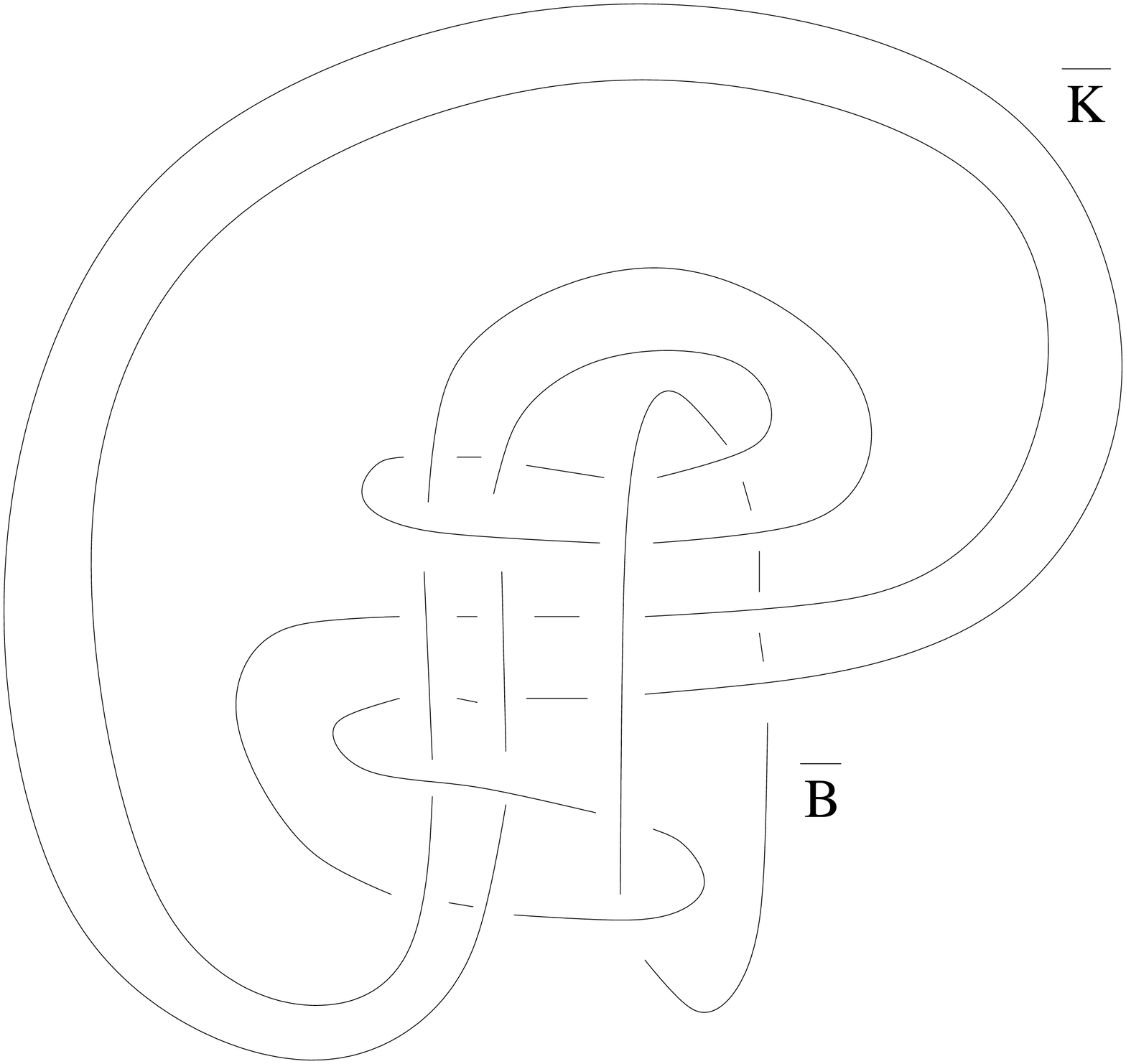}$$

\end{diagram}

It is clear from Figure \ref{ribbon2}
that in the $q$-fold cyclic cover of $S^3$ branched over
$\overline B$ the knot $\overline K$ will lift to a
$q$-equivariant ribbon knot.
\end{example}

Our proof of the realization  half of Theorem \ref{eqribbonthm} is related 
to 
 Section 2 of \cite{DL} where a construction was given to obtain a periodic knot
with
 desired Alexander polynomial.

\begin{theorem}  Given a polynomial $a(g,t) \in \ZqZ$ with $a(g,1) = 1$, 
there is a $q$-equivariant ribbon knot whose Murasugi polynomial is  
$a(g,t)a(g^{-1},t^{-1})$.
\end{theorem}

\begin{proof}  By Theorem \ref{Murasugicompute}, we must find a 
4-component link $S_1 \cup S_2 \cup \overline B \cup
\overline K$ satisfying the conditions of Remark  \ref{conditions} so that
$\lk_{\Z/q \times \Z}(S_1,S_2) = a(g,t)$.  

Let $a (g,t) =  1 + \Sigma_{i = 0}^{q-1} ~ h_i(t)g^i$, where
$h_i(t) = \Sigma_j~ a_{i,j}t^j$. 
Since $a(g,1) = 1$, 
$ \Sigma_j a_{i,j} = 0$ for all $i $.
In Example \ref{exequiv} we had $a(g,t) = 1 + h_0 (t) + g h_1(t), \ h_0(t) =
0$, 
and $h_1(t) = -1 + t$.

The desired link is indicated by the following two figures.  The polynomial
$h_i(t)$ determines
what goes in Box $i$ below.

\begin{diagram}
\label{quotientlink}
$$
 \epsfysize 2in 
\epsffile{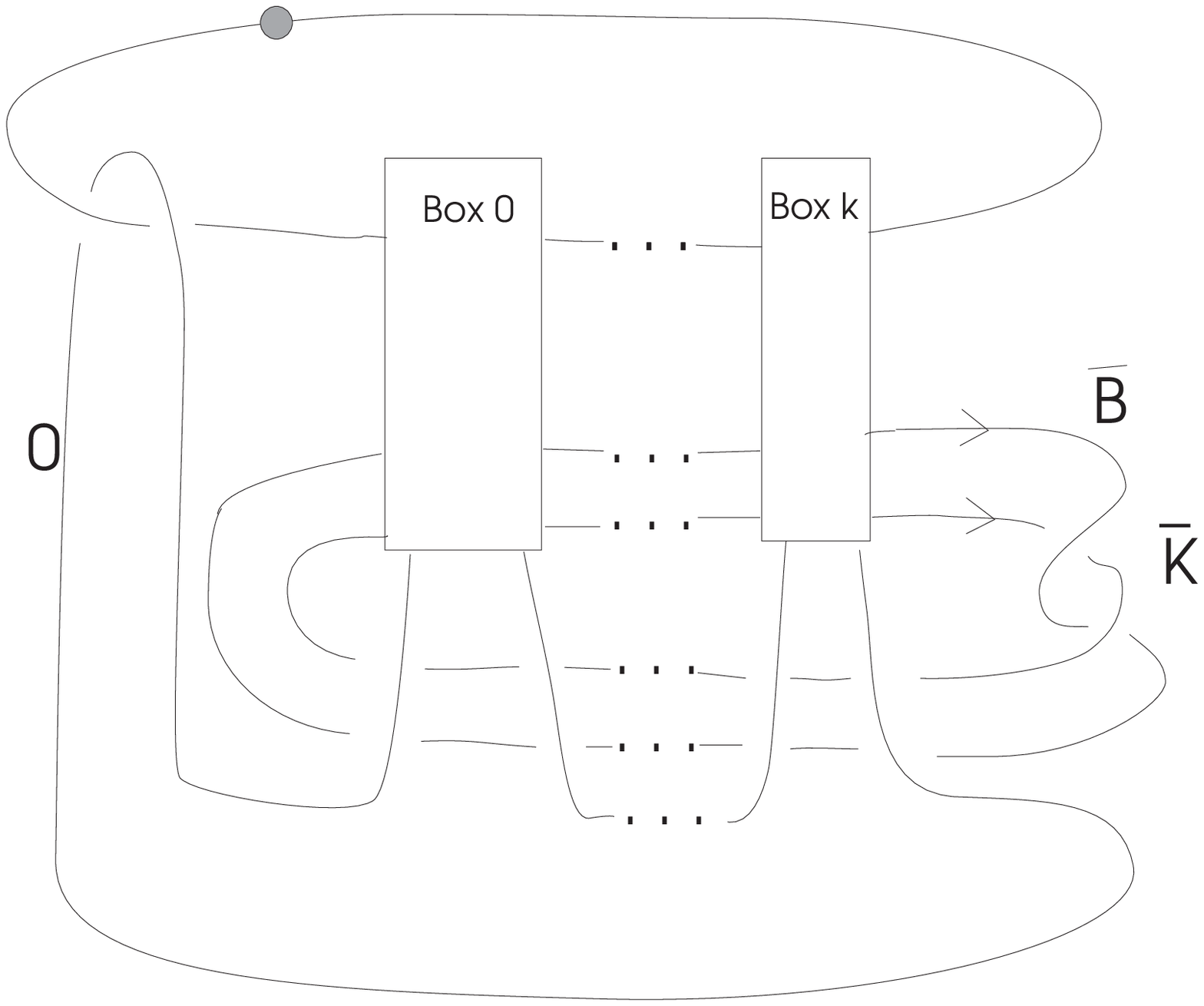}$$ 
\end{diagram}

Figure \ref{box2} illustrates Box 2 where $h_2(t) = 3t^4 - t^2 -
2 t ^{-1}$.  

\begin{diagram}
\label{box2}

$$
 \epsfysize 2in 
\epsffile{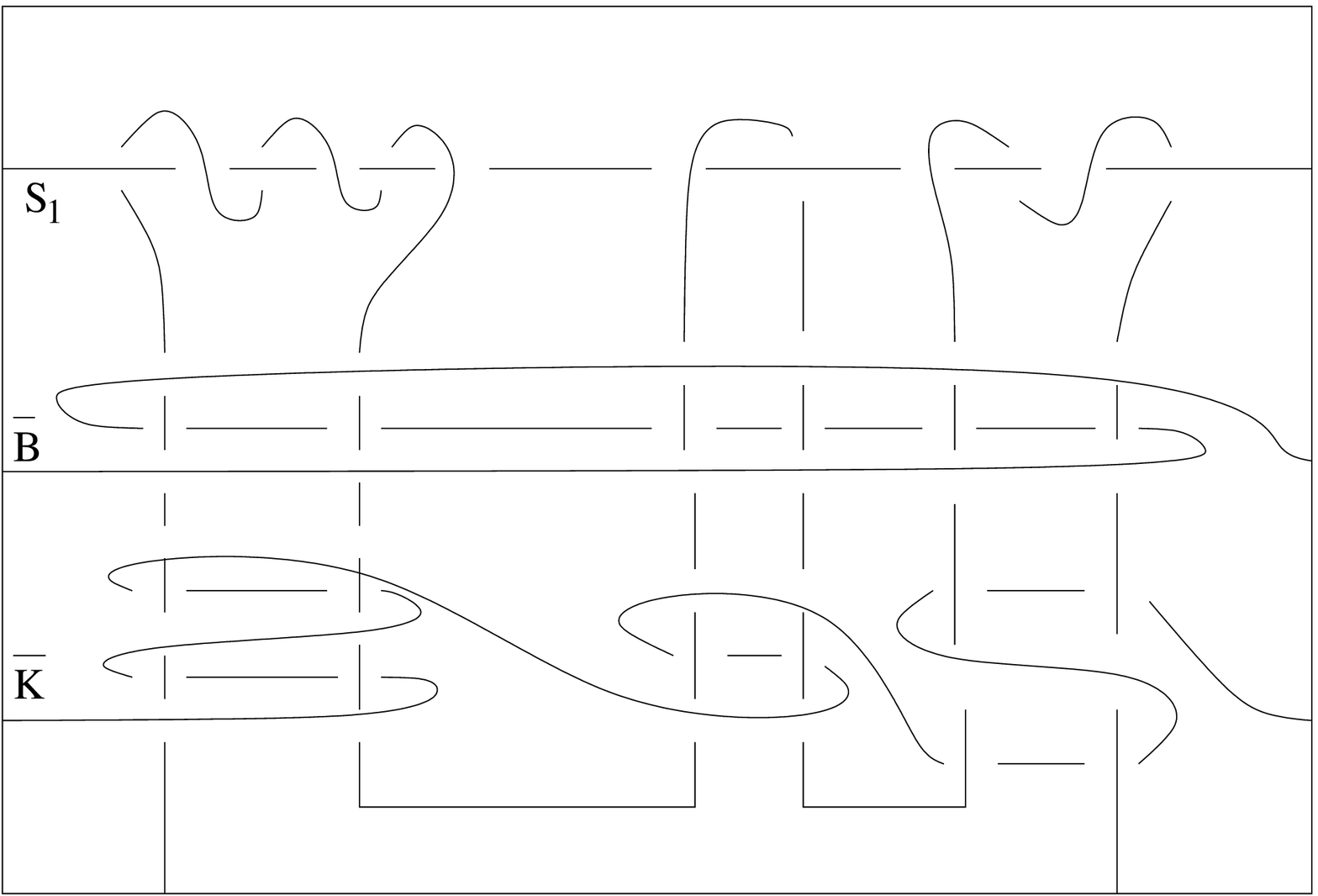}$$ 
\end{diagram}

Here $a_{2,4} = 3$ determines the shape of the the ``cactus'' on the left;
the 2 determines the linking with $\overline
B$, the 4 determines the linking with $\overline K$ and the 3 determines the
linking with $S_1$.  With this
construction $\lk_{\Z_q
\times \Z}(S_1,S_2) = a(g,t)$.  Since $ \Sigma_j a_{i,j} = 0$ for all $i$,
$\overline B$ bounds a disk (below the
boxes) in the complement of $S_1 \cup S_2$.  The other conditions of
Remark \ref{conditions} are clearly satisfied.  

The desired $q$-equivariant ribbon knot is given by adding a 1-handle along
$S_1$ and 2-handle along $S_2$, and then taking the
$q$-fold branched cover branched along $\overline B$.  

\end{proof}

\section{Background on Reidemeister torsion} \label{Background on
Reidemeister torsion}

Call a chain complex $C$ over a ring $A$ {\em finite} if
$\bigoplus C_i$ is a finitely generated module, {\em free} if
each $C_i$ is a free module, {\em based} if each $C_i$ is a free
module equipped with a basis, and {\em acyclic} if each
$H_i(C)$ is zero.  For a finite, based, acyclic complex $C$ over a ring
$A$, one can define\footnote{The definition of torsion is given,
for example,
 in \cite{Milnor1966}.
Torsion satisfies the following three axioms:
\begin{enumerate}
\item  If $C =\{\dots \to 0 \to C_1 \xrightarrow{\partial_1} C_0 \to 0 \to
\dots \}$ and $M$ is the matrix
of $\partial_1$, then $\tau(C) = [M]$.
\item $\tau(C \oplus C') = \tau(C)  \tau(C')$
\item If $f : C \to C'$ is a chain isomorphism, then
$$ \tau(C') 
\prod_{i \text{ odd}} 
[f_i : C_i \to C_{i'}]  = \tau(C) 
\prod_{i \text{ even}} 
[f_i : C_i \to
C_{i'}] 
$$
\end{enumerate}
} the {\em
torsion}
$\tau(C)
\in
\widetilde{K}_1(A)$,  where $\widetilde{K}_1(A) = K_1(A)/\{\pm I\}$,
$K_1(A) =
GL(A)/[GL(A),GL(A)]$,
and $GL(A) = \cup_{n=1}^{\infty} GL_n(A)$, 
where $GL_n(A) \hookrightarrow GL_{n+1}(A)$ via $M \mapsto 
\begin{pmatrix}
M & 0 \\
0 & 1
\end{pmatrix}$.
When
$A$ is a commutative ring there is a determinant homomorphism 
$$
K_1(A) \longrightarrow A^{\times} = \{a
\in A : ab = 1 \text{ for some } b \in A\},
$$
and so in the commutative case the torsion defines an element in
the units of $A$ modulo $\{\pm 1\}$.  In this paper, we only use the case of
$A$ commutative and we only
use the image of $\tau(C)$ in $A^\times$, so the reader who is queasy
about algebraic $K$-theory may
replace $K_1(A)$ by $A^\times$ if so desired.

A basic fact about torsion is that if 
$$
0 \to C' \to C \to C'' \to 0
$$ 
is a short exact sequence of finite, based, acyclic chain complexes, and if for
each $i$, the basis of $C_i$ is the union
of the basis for $C'_i$ and a set whose image is the basis of $C''_i$, then 
$$
\tau(C) = \tau(C')\cdot \tau(C'')
$$

When $R \to A$ is a ring homomorphism and $C$ is a finite, free chain
complex over $R$ so that 
$A \otimes_R C$ is acyclic, one can define the {\em algebraic Reidemeister
torsion} 
$$\Delta(C ;A) \in \widetilde K_1(A)/\widetilde K_1(R)
$$
by assigning any $R$-basis to $C$ and defining $\Delta(C ;A)$ to be the image
of $\tau(A \otimes_R C)$.  Actually,
in this paper, $C$ will be equipped with an equivalence class of bases, where
the torsion of the change of
basis sits in a subgroup $S \subset \widetilde K_1(R)$, in which case
$\Delta(C ;A) \in\widetilde  K_1(A)/S$ is
well defined.

We shall be concerned with the following geometric situation.  Let $X$
be a finite $CW$-complex equipped with a map $X \to BG$, where $BG$ is a
connected CW complex with $\pi_1(BG) = G$ and
all higher homotopy groups zero.  Let
$\widetilde{X}$ be the corresponding $G$-cover of $X$.  Then
$C(\widetilde{X})$ is a
finite, free $\Z G$-chain complex, and by choosing a cell in
$\widetilde{X}$ above each cell in 
$X$ and choosing an orientation for each such cell, $C(\widetilde{X})$
becomes a finite, based chain complex over $\Z G$.  
Suppose $\Z G \longrightarrow A$ is a map of
rings and that $C(X; A) = A \otimes_{\Z G} C(\widetilde{X})$ 
is acyclic.  In this case one defines the {\em Reidemeister torsion}
$\Delta(X ; A) \in K_1(A) /  \pm \! G$, where $\pm  G$
refers to the subgroup generated by the image of the elements $\pm g \in
GL_1(\Z G)$ where $g \in G$. \par
Let $f : M \to X$ be a cellular map of $CW$-complexes.  Let $C(f ; A)$ be the
algebraic mapping cone of $f_\# : C(M ; A) \to C(X ; A)$.  Here $C_i(f;A) =
C_{i-
1}(M;A) \oplus
C_{i}(X;A)$ and
$$\partial = 
\begin{pmatrix}
-\partial_M & 0 \\
 f_{\#} & \partial_X
\end{pmatrix}.
$$
There is a short exact sequence of chain complexes
$$0 \to C_{*}(X;A) \to C_*(f;A) \to C_{*-1}(M;A) \to 0
$$
and a corresponding long exact sequence on homology 
$$ \dots \to H_{i+1}(f; A) \to H_i(M; A)\to H_i(X; A)\to H_{i}(f; A)\to
\cdots .
$$
If $C(f ; A)$ is acyclic, its torsion is denoted $\Delta(f; A)$.  \par

We need three background lemmas, which correspond to three properties of
the Alexander
polynomial: it is a polynomial, it is symmetric, and it augments to 1.

\begin{lemma}  \label{integral} Let $C$ be a finite, based chain complex over
$R$ which has
the chain homotopy type of a finite, free complex $D = \{D_{i+1} \to D_{i}\}$. 
If $R \to A$ is
a ring homomorphism so that $A
\otimes_R C$ is acyclic, then $\tau(A\otimes_R C)^{(-1)^i}$ can be
represented by a
matrix whose coefficients are in the image of $R$.
\end{lemma}

\begin{proof}  Give $D$ an $R$-basis.  Let $h: C \to D$ be a chain homotopy
equivalence and give the algebraic
mapping cone $C_i(h) = C_{i-1} \oplus D_{i}$ an $R$-basis induced from that
of $C$ and $D$.  Note that $C(h)$ is acyclic, so
$\tau(C(h)) \in \widetilde K_1(R)$.  Finally,
$$
\tau(A \otimes_R C) = \text{im }\tau(C(h))\cdot\det(A \otimes_R D_{i+1}
\to A \otimes_R D_i)^{(-1)^i},
$$
by the ``basic fact''.
\end{proof}

A {\em ring with involution} is a ring $A$
equipped with a map $- : A \longrightarrow A$ satisfying
$\overline{a + b} = \overline{a} + \overline{b}$, 
$\overline{ab} = \overline{b} \overline{a}$, and $\overline{\overline{a}} = 
a$.  An example of such is a group ring $\Z G$, equipped with
the involution  $\overline{\Sigma a_g g} = \Sigma a_g g^{-1}$.  Note
that an involution on ring $A$ induces a homomorphism  $- : K_1(A)
\to K_1(A)$.

Recall that a triad $(X; X_+, X_-)$ of dimension $n$ is a compact manifold
$X$ of dimension $n$, together with a decomposition of its boundary
$\partial X = X_+ \cup X_-$, where $X_+$ and $X_-$ are compact 
$(n-1)$-dimensional manifolds with $\partial X_+ = X_+
\cap X_- = \partial X_-$. We allow the possibility of $X_+$ and/or $X_-$ being
empty.  A map of triads
$f : (M; M_+, M_-) \to (X; X_+, X_-)$ is a map $f : M \to X$   which restricts
to  maps
$M_+ \to X_+$ and $M_- \to X_-$.

\begin{lemma} \label{torsion duality}  Let $X$ be a compact, oriented
manifold of dimension
$n$.  Let $X \to BG$ be a continuous map and $\Z G \to A$ a map of 
rings with involution.
\begin{enumerate}

\item  If $X$ is closed and $C(X ; A)$ is acyclic, then
$\Delta(X ; A) = \overline{\Delta(X ; A)}^{(-1)^{n-1}}$.
\item  If $C(X ; A)$ and $C(\partial X ; A)$ are acyclic, 
then  $$\Delta(\partial X; A) = \Delta(X ; A) \overline{\Delta(X ; A)}^{(-
1)^n}.$$
\item  If $f : M \to X$ is a proper, degree one map between compact 
manifolds which restricts to a homotopy equivalence
$\partial M \to \partial X$ and if $C(f ; A)$ is acyclic, then
$\Delta(f ; A) = \overline{\Delta(f ; A)}^{(-1)^{n-1}}$.
\item  If  $f : (M; M_+, M_-) \to (X; X_+, X_-)$ is a degree one map of
triads, and if $f |_{M_-} : M_- \to X_-$ is a homotopy equivalence, then
$\Delta(f|_{M_+}; A) = \Delta(f ; A) \overline{\Delta(f ; A)}^{(-1)^n}$. 
\end{enumerate}
\end{lemma}

\begin{proof}
(1)  The proof follows Milnor \cite{Milnor1962} and Wall
\cite[Chapter 2]{Wall}.
Define $C^*(X ;A) = \text{Hom}_A(C(X) , A)$, and give it the dual basis.
Then Poincar\'e duality gives a chain equivalence
$C^*(X ; A) \to C(X; A)$ whose mapping cone has trivial torsion (this is
by the dual cell proof of Poincar\'e duality).  It follows that the domain
and range have equal torsion, and after taking into account the dual maps
and dimension shift the torsion of the domain is 
$\overline{\Delta(X ; A)}^{(-1)^{n-1}}$.  (2) follows from a relative
version of (1).  Finally (3) and (4) follow because 
$C(f ; A)$ satisfies Poincar\'e duality since $f$ has degree one. 
\end{proof}

For a group $G$, let $\epsilon : \Z G \to \Z$ be the augmentation map
$\epsilon (\Sigma~ a_g g) =
\Sigma~ a_g $ and 
$S_G = \{\alpha
\in
\Z G  : \epsilon(\alpha) = 1
\}$.

The following lemma is proven in \cite{SmithJR}.

\begin{lemma}  \label{smith} Let $C$ be a finite, free chain complex over
$\Z G $ for a
finitely generated abelian group $G$.  Then 
$\Z \otimes_{\Z G } C$ is acyclic if and only if $S_G^{-1}C = S_G^{-1}\Z G 
 \otimes_{\Z G } C$ is acyclic.
\end{lemma}

The localization $S_G^{-1}\Z G $ is defined abstractly, but in cases of
interest there are
explicit realizations:
$$S^{-1}_\Z\Z[\Z] = \{f(t)/h(t) \in \Q(t) : f,h \in \Z[t,t^{-1}], h(1) = 1\}$$
$$S^{-1}_{\Z/q \times \Z}\Z[\Z/q \times\Z] = \{f(g,t)/h(g,t) \in
\Q[\Z/q]\times \Q(t) : f,h \in
\Z[\Z/q][t,t^{-1}], h(1,1) = 1\}$$
$$S^{-1}_{\Z \times \Z}\Z[\Z \times\Z] = \{f(g,t)/h(g,t) \in \Q(g,t) : f,h \in
\Z[g,g^{-1},t,t^{-1}], h(1,1) = 1\}$$

\begin{lemma}  \label{Alexpolytorsion} Let $K$ be a knot. Let $L$ be a 
2-component link with 
linking number one.    
\begin{enumerate}
\item There are tubular neighborhoods $N(K)$ and $N(L)$ and proper maps
$f: S^3 - N(K) \to S^1 \times D^2 $ and
$h :S^3 -N(L) \to S^1 \times S^1 \times [-1,1]$ inducing isomorphisms on
homology and homeomorphisms on the respective
boundaries.
\item The chain complexes $C(f; S^{-1}_\Z\Z[\Z])$ and $C(h; S^{-1}_{\Z
\times \Z}\Z[\Z \times\Z])$
are acyclic,   
$\Delta(f; S^{-1}_\Z\Z[\Z]) = \Delta_K(t)$, and $\Delta(h; S^{-1}_{\Z
\times \Z}\Z[\Z \times\Z]) =
\Delta_L(t)$.
\item If $L = \overline{B} \cup \overline{K}$ is a link with $\overline{B}$
unknotted, then $\Delta(h; S^{-1}_{ e \times \Z }\Z[\Z/q \times\Z])$
represents the Murasugi polynomial
of the period $q$ knot given by taking the $q$-fold cover of $S^3$ branched
over $\overline{B}$. 

\end{enumerate}
\end{lemma}

\begin{proof}  1.  Let
$N(K)$ be the open disk bundle of a tubular neighborhood of
$K$.  There is a commutative diagram
$$\begin{CD}
H^1\partial(S^3 - N(K)) @<<< H^1 (S^1 \times S^1) \\
@AAA @AAA \\
H^1 (S^3 - N(K)) @<<< H^1 (S^1 \times D^2)
\end{CD}
$$
where the vertical maps are induced by inclusions of spaces, and the
horizontal maps are isomorphisms given by Alexander duality.
Realize the top horizontal map by a homeomorphism, and extend this to a
map $f$ using that $S^1
\times D^2$ is a $K(\pi,1)$-space (see \cite[V.6.9]{Wh}). The construction
of $h$ is similar, and linking number one is used to
see that $h$ gives an isomorphism on homology.

2.  The acyclicity follows from Lemma \ref{smith}.  The relationship between
Alexander polynomials and
Reidemeister torsion is well-known (see
\cite{Milnor1962} and \cite{Tu}).  In our case it follows since $C(f;\Z[\Z])$
has the chain
homotopy type of a presentation of $H_1$ of the infinite cyclic cover, and
the determinant of a
presentation is the Alexander polynomial.  The argument for the Alexander
polynomial of the link
is similar. 

3.  In the statement
$$
S_{e\times \Z} = \{\alpha(g,t) \in\Z[\Z \times \Z] : \alpha(g,1) =1\}
\subset S_{\Z \times \Z}.
$$
We see $\Z \otimes_{\Z[\Z]} C(h;\Z[\Z/q \times\Z])$ is acyclic since
$\overline{B}$ is unknotted.  Hence by Lemma \ref{smith}, $
S^{-1}_{ e \times \Z }C(h;\Z[\Z/q \times\Z])$ is acyclic.  This result then
follows from Part 2., since the Murasugi polynomial
is represented by the Alexander polynomial of the link.
\end{proof}

\section{Proofs of Theorems \ref{eqslicethm} and \ref{eqribbonthm}}
\label{Application to equivariant slice/ribbon knots}

Let $K$ be a knot of period $q$ which is equivariant slice.  Let $D_K$ be the
slice disk, invariant under the group action, and let $D_B$ be the 
fixed set of the action on $D^4$.  Let $B = D_B \cap S^3$.  Then the linking
number of $K$ and
$B$ in $S^3$ is one \cite{N}.  (The basic idea is that the linking number
equals the
intersection number of $D_K$ and $D_B$.  Intersection points correspond to
fixed points in $D_K$, and there can be only one.)  Let
$\overline{B}^4=B^4/(\Z/q)$; this a 
topological ball with boundary $S^3$.  Let $D_{\overline{B}} \subset
\overline{B}^4$  denote the image
of the fixed set in the quotient,  $D_{\overline{K}}\subset \overline{B}^4$
the quotient of the
slice disk, and $\overline{B} =D_{\overline B} ~\cap ~S^3 $.  
A Mayer-Vietoris argument shows that
$H_1(\overline{B}^4 ~-~ D_{\overline{B}}~ -~ D_{\overline{K}}) =
\Z \times \Z$ (let $\overline B^4 = ( N(D_{\overline{B}}) \cup
N(D_{\overline{K}})) \cup (\overline{B}^4 ~-~
D_{\overline{B}}~ -~ D_{\overline{K}})$).  Let $\partial_+ D^2$ and
$\partial_- D^2$ denote the upper and lower
semicircles of $S^1 \subset D^2$.  Then there are tubular neighborhoods
and a map of triads
$$f :  (\overline{B}^4 - N(D_{\overline{B}}) - N(D_{\overline{K}});~ 
S^3 - N(\overline{B}) - N(\overline{K}),~
\partial (\overline{N(D_{\overline{B}}) \cup N(D_{\overline{K}}))})$$
$$\to
(D^2 \times S^1 \times S^1 ; ~\partial_+D^2 \times S^1 \times S^1, 
\partial_-D^2 \times S^1 \times S^1),$$
inducing a homeomorphism on the third component of the triad.

  The proof of Theorem \ref{eqslicethm} follows from  Lemma
\ref{torsion
duality}(4) (with $A = S^{-1}_{e \times \Z} \Z[\Z/q\times \Z]$), 
Lemma \ref{Alexpolytorsion}(3), and clearing denominators.  

Let $K$ be a $q$-equivariant ribbon knot with Alexander polynomial
$\Delta (t)$  and Murasugi polynomial $\Delta_{\Z /q}(g,t)$. Since $K$ is
equivariant ribbon,
$C_*(f;\ZqZ)$ is a finite, free $\ZqZ$-chain complex with at most one
nonzero homology group $M
= H_2(f; \ZqZ)$.  We claim any such chain complex is chain homotopy
equivalent to a finite, free
chain complex of the form $D_3 \to D_2$.  The first step is to use that a
surjective map to a
free module splits, so that $C_*(f;\ZqZ)$ is chain homotopy equivalent to a
finite, free
complex $C'_*$ which is zero below degree 2.  Hence $M$ has finite
homological dimension, and is
$\Z$-free by Poincar\'e duality, and hence by Rim's theorem \cite[Section
4]{Rim} is projective
over $\Zq$.  Then we see $M$ is homological dimension 1, since
$$0 \to \Z[t,t^{-1}] \otimes_\Z M  \xrightarrow{\cdot(1-t) \otimes id }
\Z[t,t^{-1}] \otimes_\Z M
\to \Z \otimes_\Z M \to 0
$$
is a projective $\ZqZ$-resolution of length 1. If necessary, we may sum on a
projective complement
to
$\Z[t,t^{-1}] \otimes_\Z M  \xrightarrow{\cdot(1-t) \otimes id } \Z[t,t^{-
1}] \otimes_\Z M$ to
find an exact sequence
$$ 0 \to D_3 \to D_2 \to M\to 0
$$
where $D_2$ and $D_3$ are finitely generated free over $\ZqZ$.   The
fundamental lemma of
homological algebra then says that $C'_*$ and $D_3 \to D_2$ are chain
homotopy equivalent.

 The existence of a polynomial 
$a(g,t)$ satisfying the requirements of Theorem \ref{eqribbonthm} then
follows  from 
Lemma \ref{torsion duality}(4),  Lemma \ref{Alexpolytorsion}(3), and
Lemma \ref{integral}.

Finally we note that we can use Poincar\'e duality to give an alternate
computation
of the Murasugi polynomial of $q$-equivariant ribbon knot that we
constructed in Section
\ref{Construction of
equivariant ribbon knots}. 

\medskip

{\em Acknowledgement:} The authors would like to thank Charles Livingston for useful conversations throughout this project.



\end{document}